# A note on the sign (unit root) ambiguities of Gauss sums in index 2 and 4 cases

YANG Jing[*]

Department of Mathematical Science, Tsinghua University, Beijing, 100084, China

Email: jingyang@math.tsinghua.edu.cn

XIA Lingli

Basic Courses Department of Beijing Union University, Beijing, 100101, China

Email: lingli@buu.edu.cn

**Abstract**

Recently, the explicit evaluation of Gauss sums in the index 2 and 4 cases have been given in several papers (see [2, 3, 7, 8]). In the course of evaluation, the sigh (or unit root) ambiguities are unavoidably occurred. This paper presents another method, different from [7] and [8], to determine the sigh (unit root) ambiguities of Gauss sums in the index 2 case, as well as the ones with odd order in the non-cyclic index 4 case. And we note that the method in this paper are more succinct and effective than [8] and [7].

**Keywords:** Gauss sum, Teichmüller maps and characters, Stickelberger's congruence, Stickelberger's theorem.

**MSC:** 11L05, 11T24

## 1 Introduction

Gauss sums are one of the most important and fundamental objects and tools in number theory and arithmetical geometry. The explicit evaluation of Gauss sums is a very important but very difficult problem. Let $p$ be a prime number, and let $N$ ($\geqslant 2$) be a natural number such that $(N, p) = 1$. Take $f = \mathrm{ord}_N(p)$ be the multiplicative order of $p$ in $(\mathbb{Z}/N\mathbb{Z})^*$, i.e. $f$ is the smallest positive integer such that $p^f \equiv 1 \pmod{N}$. For a multiplicative character $\chi$ with fixed order $N$ over the finite field $\mathbb{F}_{p^f} = \mathbb{F}_q$, let $\chi(0) = 0$. Thus, for $1 \leqslant r \leqslant N-1, 1 \leqslant v \leqslant p-1$, the Gauss sums $G(\chi^r, v)$ over $\mathbb{F}_q$ related to $\chi^r$ and $v$ are defined by

$$G(\chi^r, v) = \sum_{x \in \mathbb{F}_q} \chi^r(x) \zeta_p^{vT(x)}. \tag{1.1}$$

where $\zeta_p = e^{2\pi i/p}$, $i = \sqrt{-1}$, and $T$ denotes the trace map from $\mathbb{F}_q$ onto $\mathbb{F}_p$.

Since $G(\chi^r, v) = \overline{\chi^r(v)} G(\chi^r, 1)$ we can concentrate on $G(\chi^r, 1)$, which we briefly denote by $G(\chi^r)$ in the following text for simplicity. And then $G(\chi) \in \mathbb{Z}[\zeta_{Np}]$, which is the integral ring of the cyclotomic field. As we known, the Galois group $\mathrm{Gal}\left(\mathbb{Q}(\zeta_N, \zeta_p)/\mathbb{Q}\right)$ is canonically isomorphic to group $(\mathbb{Z}/Np\mathbb{Z})^* \cong (\mathbb{Z}/N\mathbb{Z})^* \times (\mathbb{Z}/p\mathbb{Z})^*$.

In 1890, L. Stickelberger gave a series profound results of Gauss sums in paper [5], which made far-reaching influence for later studies of Gauss sums. One of these famous results is called "Stickelberger's congruence" by later people.

---

[*]Corresponding author: YANG Jing, Email: jingyang@math.tsinghua.edu.cn

Recently, a series papers [2, 7, 3, 8] have studied the problem of explicit evaluation of Gauss sums in the index 2 and 4 case, i.e., the case that $[(\mathbb{Z}/N\mathbb{Z})^* :< p >] = \varphi(N)/f = 2, 4$, where $\varphi$ is Euler's function. Throughout the proof of these papers, one can find that there unavoidably exists the sigh (or unit root) ambiguities in the course of evaluation. The reasons for the occurrence of sigh (unit root) ambiguities lie that: (1) quadric equations always have solutions with "$\pm$" ambiguities; (2) when one finds generator element for a principal ideal in the ring, there is inevitably an ambiguity of the unit, exactly it's an ambiguity of the unit root in the ring. The methods to solve the problem of sigh ambiguity in the above papers are that: lifting both sides of the formulas for Gauss sums to the power of $l_1^{r_1}$ ($l_1 \mid N$ is odd), and then taking both sides modulo $l_1$. Thus, the sigh ambiguities can be distinguished. However, when $l_1 = 1, 3$, $\sqrt{-l_1} \in K$, which is the decomposition subfield of $p$, the integral ring $O_{\mathbb{Q}(\sqrt{-l})} \subset O_K$ has four or six unit roots, i.e., $<i>$ or $<\zeta_6>$, they can not be distinguished by the congruences modulo $l_1$. So, we need another method to determine the ambiguities in these cases.

In this paper, using Stickelberger's congruence, another method to determine the sigh (unit root) ambiguities of Gauss sums different from [7, 8] is presented. First, some preliminary results are given in Section 2, including fundamental knowledge of Teichmüller maps and characters and the representation of Stickelberger's congruence in detail. Then, in Section 3, we give the main theorem to solve the problem of sigh (unit root) ambiguities of Gauss sums in the index 2 case, as well as some concrete examples. Finally, in Section 4, for the Gauss sums with odd order $N$ in the non-cyclic index 4 case, we do the same things as the ones in Section 3.

## 2 Preliminary Results

### 2.1 Teichmüller Maps and Characters

First, we introduce the Teichmüller maps and characters, and then determine the character $\chi$, related with the Gauss sums, of order $N$. Let $L = \mathbb{Q}(\zeta_{q-1})$, and take $\mathfrak{P}$ be one of prime ideal factors of $p$ in $O_L$. By the extension theories of cyclotomic fields, we know that

$$O_L/\mathfrak{P} \cong \mathbb{F}_q = \mathbb{F}_{p^f}. \tag{2.1}$$

For $0 \leqslant j \leqslant q-1$, since $q = \prod_{j=1}^{q-1}(1 - \zeta_{q-1}^j)$, $(q-1, p) = 1$, thus $1 - \zeta_{q-1}^j \notin \mathfrak{P}$, i.e. $\zeta_{q-1}^i \notin \zeta_{q-1}^j + \mathfrak{P}$ ($i \neq j$). So $\zeta_{q-1}^j$ respectively belong to the different cosets (mod $\mathfrak{P}$) when $j$ runs from 1 to $q-1$.

Now define Teichmüller character $\omega_\mathfrak{P}$ and a character $\chi$ with order $N$ as follow:

$$\begin{aligned} \omega_\mathfrak{P}(t) &= \zeta_{q-1}^j \quad (\forall t \notin \mathfrak{P}) \\ \chi(t) &= \chi_\mathfrak{P}(t) = \omega_\mathfrak{P}^{\frac{q-1}{N}}(t) \end{aligned} \tag{2.2}$$

where $\zeta_{q-1}^{j(q-1)/N}$ denote the unique $(q-1)$-th unit root in $\bar{t}$. Moreover, by the above correspondences, we can lift the isomorphism (2.1) to the following maps, which are called Teichmüller maps

$$\begin{aligned} \phi: \quad \mathbb{Z}[\zeta_{q-1}] &\to \mathbb{F}_q = <\gamma> \cup \{0\} \\ \sum_{i=0}^{q-1} a_i \zeta_{q-1}^i &\mapsto \sum_{i=0}^{q-1} b_i \gamma^i \end{aligned} \tag{2.3}$$

where $a_i \equiv b_i \pmod{p}$, $b_i \in \mathbb{F}_p$, and the selection of $\gamma$, the primitive element of $\mathbb{F}_q^*$, is determined by $\gamma \equiv \zeta_{q-1}$ (mod $\mathfrak{P}$). It is easy to know that $\phi$ is a surjective homomorphism. Particularly, if $K$ is the fixed subfield of isomorphism $\sigma_p$, i.e. for $x \in K$ one has $\sigma_p(x) = x$, then by Teichmüller map $\phi$, $\phi(x)^p = \phi(x)$, thus $\phi(x) \in \mathbb{F}_p$.

From now on, we always view $\chi$ as $\chi_\mathfrak{P}$ defined in (2.2). Next, we consider the Gauss sum $G(\chi)$ related to $\chi = \chi_\mathfrak{P}$.

## 2.2 Stickelberger's congruence

In 1890, L. Stickelberger stated a series profound results of Gauss sums in [5], which made long influence for the research after. A famous result among them, called "The Stickelberger's Theorem". It gives a revealment of the factorization for the Gauss sums in the integral ring of cyclotomic field. As we known, for a multiplicative character $\chi$ with order $N$ over $\mathbb{F}_q$, $G(\chi)^N \in \mathbb{Z}[\zeta_N] = O_M$, where $M = \mathbb{Q}(\zeta_N)$. And the Galois group

$$\mathrm{Gal}(M/\mathbb{Q}) = \{\sigma_a | 1 \leqslant a \leqslant N-1, (a, N) = 1\}$$

is canonically isomorphic to multiplicative group $(\mathbb{Z}/N\mathbb{Z})^*$.

**Theorem 2.1** (Stickelberger's theorem, see in [1, §11.2] or [4, §11.3]). *Let $G(\chi)$ be the Gauss sum with order $N$ over $\mathbb{F}_q$, $M = \mathbb{Q}(\zeta_N)$. Then, there exists a prime ideal factor $P$ of $p$ in $O_M$ such that*

$$G(\chi)^N O_M = P^\theta,$$

*where $\theta = \sum\limits_{\substack{a=1 \\ (a,N)=1}}^{N-1} a\sigma_a^{-1} \in \mathbb{Z}[G]$ is usually called :"Stickelberger Element".*

Actually, in [5] L. Stickelberger also gave a more exact congruence relation for Gauss sums in cyclotomic field $\mathbb{Q}(\zeta_N, \zeta_p)$, which is called "Stickelberger's congruence", and Theorem 2.1 is just the corollary of it.

For $0 \leqslant a < q-1$ ($q = p^f$), let the $p$-adic expression of $a$ be

$$a = c_0 + c_1 p + \cdots + c_{f-1} p^{f-1}$$

where $c_i \in \{0, 1, \cdots, p-1\}$. And we define that

$$S(a) = c_0 + \cdots + c_{f-1} \in \mathbb{Z}, \quad t(a) = c_0! c_1! \cdots c_{f-1}! \in \mathbb{Z}.$$

And for any $b \in \mathbb{Z}$, assume that $S(b) = S(b')$, $t(b) = t(b')$, where $1 \leqslant b' < q-1$, $b' \equiv b \pmod{q-1}$.

**Theorem 2.2** (Stickelberger's congruence, see in [1, §11.2]). *Let $\chi_{\mathfrak{P}}$ define as (2.2) and $v_{\mathfrak{P}}$ denote the the exponential valuation for prime divisor $\mathfrak{P}$. Then, for $1 \leqslant a \leqslant N-1$, $n = \frac{q-1}{N}$,*

$$v_{\mathfrak{P}}(G(\chi^{-a})) = S(an).$$

*More exactly, let $\pi = \zeta_p - 1$, $(v_{\mathfrak{P}}(\pi)) = 1$, then*

$$G(\chi^{-a}) \equiv -\frac{\pi^{S(an)}}{t(an)} \pmod{\mathfrak{P}^{S(an)+1}}$$

we remark that the character $\chi_{\mathfrak{P}}$ in Theorem 2.2 is determined by a properly selected $\mathfrak{P}$, which is a prime ideal factor of $p$ in $\mathbb{Q}(\zeta_{q-1})$. The more details for "How to select $\mathfrak{P}$?" will be given in Section 3 and 4.

With the preparation above, we next consider the determination of the sigh (unit root) ambiguities of Gauss sums by Stickelberger's congruence, which is a different method from [7, 8]. And we always assume that $p \neq 2$ except so far as otherwise expressly stated.

# 3 Sigh (unit root) ambiguities of Gauss sums in index 2

In [8], explicit evaluation of Gauss sums in index 2 case has been given. And we find that there are the following four forms for Gauss sums in index 2 case:

$$\begin{aligned}&(1). \ G(\chi) = \varepsilon p^{\frac{f}{2}}; \\ &(2). \ G(\chi) = \varepsilon\sqrt{p^*} \cdot p^{\frac{f-1}{2}} \quad (p^* = (-1)^{\frac{p-1}{2}} p); \\ &(3). \ G(\chi) = 2^{-y} p^x (a + b\sqrt{-L}) \quad (x, y, L \in \mathbb{Z}_+ \cup \{0\}); \\ &(4). \ G(\chi) = 2^{-y} \sqrt{p^*} \cdot p^x (a + b\sqrt{-L}).\end{aligned} \quad (3.1)$$

To determine $\varepsilon, a$, [8] considered the congruences modulo $L$ or certain odd factor of $L$. However, when $L = 1, 3$, $\mathbb{Q}(\sqrt{-L})$ has 4 or 6 unit roots, and by the congruences modulo $L$ (or its factor), one can not distinguish these unit roots. Additionally, when $L = 2$, by the congruences modulo $L$, one can not distinguish $\pm 1$ at all. Thus, in these cases, we need another method to solve the problem of sigh (unit root) ambiguities for Gauss sums in index 2 case.

First, we consider the Gauss sums as form (1) in (3.1), i.e. $G(\chi) = \varepsilon p^{\frac{f}{2}}$. By $\pi = \zeta_p - 1$ and congruence

$$p = \prod_{i=1}^{p-1}(1-\zeta_p^i) \equiv \pi^{p-1}\prod_{i=1}^{p-1} i \equiv -\pi^{p-1} \pmod{\mathfrak{P}^p}, \tag{3.2}$$

we can determine $\varepsilon$ directly by Stickelberger congruence, i.e.,

$$G(\chi) \equiv \varepsilon(-1)^{\frac{f}{2}}\pi^{(p-1)\frac{f}{2}} \equiv \frac{-\pi(-n)}{t(-n)} \pmod{\mathfrak{P}^{s(-n)+1}}.$$

By Stickelberger's theorem, $s(-n) = (p-1)\frac{f}{2}$, $(\pi) = \mathfrak{P}$, we have that

$$\varepsilon \equiv (-1)^{\frac{f}{2}+1}/t(-n) \pmod{\mathfrak{P}}. \tag{3.3}$$

Since $\varepsilon, (-1)^{\frac{f}{2}+1}, t(-n)$ are all belong to field $K$, which is the fixed subfield of $\sigma_p$ in $\mathbb{Q}(\zeta_p, \zeta_N)$, acting Teichmüller map $\phi$ to both sides of the (3.3), we have that

$$\varepsilon \equiv (-1)^{\frac{f}{2}+1}/t(-n) \pmod{p}. \tag{3.4}$$

Next, for the Gauss sums as form (2) in (3.1), i.e., $G(\chi) = \varepsilon\sqrt{p^*} \cdot p^{\frac{f-1}{2}}$, by Stickelberger's congruence,

$$\sqrt{p^*} \equiv -\pi^{\frac{p-1}{2}}/(\frac{p-1}{2})! \pmod{\mathfrak{P}^{\frac{p+1}{2}}}, \tag{3.5}$$

where $\sqrt{p^*}$ equals to the quadratic Gauss sum over $\mathbb{F}_p$. Thus,

$$G(\chi) \equiv -\varepsilon(-1)^{\frac{f-1}{2}}\pi^{\frac{p-1}{2}}\pi^{\frac{f-1}{2}(p-1)}/(\frac{p-1}{2})! \equiv -\pi^{s(-n)}/t(-n) \pmod{\mathfrak{P}^{s(-n)+1}}.$$

Similarly, let $\varepsilon_p \equiv (\frac{p-1}{2})! \pmod{p}$, then

$$\varepsilon = (-1)^{\frac{f-1}{2}} \cdot \varepsilon_p/t(-n) \pmod{p}. \tag{3.6}$$

We note that when $\varepsilon \neq \pm 1$, such as, when $\varepsilon \in <\sqrt{-1}>$, we need solve equation $x^2 + 1 = 0$ in $\mathbb{F}_p$ and properly select $\mathfrak{P}$ to determine the responding elements in $\mathbb{F}_p$ to $\phi(\sqrt{-1})$ and $\phi(-\sqrt{-1})$, and then determine $\varepsilon$.

Sequently, for the Gauss sums as form (3) in (3.1), i.e., $G(\chi) = 2^{-y}p^x(a + b\sqrt{-L})$, by Stickelberger's congruence, we have

$$(-1)^x\pi^{(p-1)x}2^{-y}(a+b\sqrt{-L}) \equiv G(\chi) \equiv \frac{-\pi^{s(-n)}}{t(-n)} \pmod{\mathfrak{P}^{s(-n)+1}}. \tag{3.7}$$

And for $\overline{G(\chi)} = \chi(-1)G(\overline{\chi}) = \chi(-1)G(\chi^{-1})$, similarly by Stickelberger's congruence, we obtain that

$$\chi(-1)(-1)^x\pi^{(p-1)x}2^{-y}(a-b\sqrt{-L}) \equiv G(\chi^{-1}) \equiv \frac{-\pi^{s(n)}}{t(n)} \pmod{\mathfrak{P}^{s(n)+1}}. \tag{3.8}$$

By Stickelberger's theorem, $b_0(p-1) = s(-n), b_1(p-1) = s(n)$, $x = \min\{b_0, b_1\}$. Let

$$T_0' = \frac{(-1)^{x+1}2^y}{t(-n)}, T_1' = \frac{(-1)^{x+1}\chi(-1)2^y}{t(n)} \quad \text{and } T_i = \begin{cases} T_i' & \text{when } x = b_i; \\ 0 & \text{otherwise.} \end{cases} (i = 0, 1) \tag{3.9}$$

Combining (3.7) and (3.8), and then simplifying, we have

$$\begin{cases} a + b\sqrt{-L} \equiv T_0 \\ a - b\sqrt{-L} \equiv T_1 \end{cases} \pmod{\mathfrak{P}}.$$

Since $\sqrt{-L}$ or $\frac{-1+\sqrt{-L}}{2}$ is the integral base of $K$, $\sqrt{-L} \in K$. Under the Action of Teichmüller map $\phi$, we have $\phi(\sqrt{-L}) \in \mathbb{F}_p$. It is not hard to get the two roots of equation $x^2 + L = 0$ in $\mathbb{F}_p$, and without loss of generality, we suppose that the two roots are $\pm\beta$. Then we can let $\phi(\sqrt{-L}) = \beta$, conditionally proper selection of $\mathfrak{P}$. So, the above congruence can be simplify as the following equation modulo $p$:

$$\begin{pmatrix} 1 & \beta \\ 1 & -\beta \end{pmatrix} \begin{pmatrix} a \\ b \end{pmatrix} \equiv \begin{pmatrix} T_0 \\ T_1 \end{pmatrix} \pmod{p}.$$

Therefore,

$$\begin{pmatrix} a \\ b \end{pmatrix} \equiv \begin{pmatrix} 1 & \beta \\ 1 & -\beta \end{pmatrix}^{-1} \cdot \begin{pmatrix} T_0 \\ T_1 \end{pmatrix} \pmod{p}. \tag{3.10}$$

Thus, we can not only determine the sigh of $a$ but also the sigh of $b$, i.e., solve the problem of sigh ambiguity for this case, conditionally proper selection of $\mathfrak{P}$.

Finally, for the Gauss sums as form (3) in (3.1), i.e., $G(\chi) = 2^{-y}\sqrt{-p^*} \cdot p^x(a + b\sqrt{-L})$, similarly, we let

$$T_0' = \frac{(-1)^x \varepsilon_p}{t(-n)}, T_1' = \frac{(-1)^x \varepsilon_p \chi(-1)}{t(n)} \quad \text{and } T_i = \begin{cases} T_i' & \text{when } x = b_i; \\ 0 & \text{otherwise.} \end{cases} \quad (i = 0, 1) \tag{3.11}$$

And let $\phi(\sqrt{-L})$ equal to $\beta$, one root of equation $x^2 + L = 0$ in $\mathbb{F}_p$. Then, $a, b$ can be determined by congruence equations (3.10) in $\mathbb{F}_p$, conditionally proper selection of $\mathfrak{P}$.

In conclusion,

**Theorem 3.1.** *For the Gauss sums in index 2 case, conditionally proper selection of $\mathfrak{P}$, the sigh (unit root) ambiguities can be solved by formula (3.4), (3.6), (3.9), (3.10) and (3.11), according to the different forms of (3.1).*

In the end of this section, we give an example, as a complement of [8], for the explicit evaluation of Gauss sums in the Case F, one of subcases of index 2 case.

**Example 1**(Case F) Let $N = 4l_1^{r_1}$, $\sqrt{-1} = i$.

(1). Firstly, take $l_1 = 5$, then $N = 4 \times 5 = 20$, $f = 4$, $h_1 = h(\mathbb{Q}(\sqrt{-5})) = 2$.

When $p \equiv 3 \pmod{4}$ and $a_1 = 1$, it belongs to Case F1 in [8]. We take $p = 3$, then, $<p> = \{1, 3, 9, 7\}$, $g_1 < p> = \{11, 13, 17, 19\}$. Solving the quadratic equations, we get

$$a^2 + 5b^2 = 3^2 \quad \Rightarrow \quad \begin{cases} a = \pm 2 \\ b = \pm 1 \end{cases}.$$

Next, we try to determine the sigh of $a, b$ by our Theorem 3.1. Since $n = (3^4 - 1)/20 = 4 = 1 \cdot 3 + 1$, $(20 - 1)n = 2 \cdot 3^3 + 2 \cdot 3^2 + \cdot 3 + 1$. So $s(-n) = 6 = b_0(p - 1)$, $s(n) = 2 = b_1(p - 1)$, $b_1 = \min\{b_0, b_1\} = 1$, $(t(n))^{-1} = (1!1!)^{-1} \equiv 1 \pmod{3}$. The two roots of equation $x^2 + 5 = 0$ in $\mathbb{F}_3$ are exactly $\pm 1$, and then we let $\phi(\sqrt{-5}) = 1$, conditionally proper selection of $\mathfrak{P}$. By (3.9),

$$\begin{pmatrix} a \\ b \end{pmatrix} \equiv \begin{pmatrix} 1 & 1 \\ 1 & -1 \end{pmatrix}^{-1} \begin{pmatrix} 0 \\ (-1)^1 \cdot 1 \cdot \chi(-1) \end{pmatrix} \equiv \begin{pmatrix} 1 \\ -1 \end{pmatrix} \pmod{3}.$$

Thus, the Gauss sum of order 20 over $\mathbb{F}_3$ and its conjugation are:

$$G(\chi) = 2(-2 - \sqrt{-5}), \; G(\overline{\chi}) = 3(-2 + \sqrt{-5}).$$

When $p \equiv 1 \pmod{4}$ and $a_1 = 1$, it belongs to Case F2(1) in [8]. We take $p = 13$, then $<p> = \{1, 13, 9, 17\}$, $g_1 <p> = \{3, 7, 11, 19\}$. So $b_0 = b_1 = 2 = \frac{f}{2}$. Then

$$G(\chi) = \varepsilon p^{\frac{f}{2}} \quad (\text{where } \varepsilon \in \{\pm 1, \pm i\}).$$

Since $n = (13^4-1)/20 = 1428 = 8\cdot 13^2 + 5\cdot 13 + 11$, $(20-1)n = 12\cdot 13^3 + (12-8)\cdot 13^2 + (12-5)\cdot 13 + (12-11)$. So $(t(-n))^{-1} = (12!4!7!1!)^{-1} \equiv 8 \pmod{13}$. And by (3.4),
$$\varepsilon \equiv (-1)^{2+1} \cdot 8 \equiv 5 \pmod{13}.$$

The two roots of equation $x^2 + 1 = 0$ in $\mathbb{F}_{13}$ are exactly $\pm 5$, and then we let $\phi(i) = 5$, conditionally proper selection of $\mathfrak{P}$. Thus, $\varepsilon = i$ and
$$G(\chi_{\mathfrak{P}}) = 13^2 i, \quad G(\overline{\chi}_{\mathfrak{P}}) = -13^2 i.$$

(2). Secondly, let $l_1 = 7$, then $N = 4 \times 7 = 28$, $f = 6$.

When $p \equiv 1 \pmod 4$ and $a_1 = 1$, it belongs to Case F2(2). We take $p = 5$, so, $<p>= \{1,5,25,13,9,17\}$, $g_1 <p>= \{3,11,17,19,23,27\}$. And $b_0 = \frac{5}{2} = \frac{f-1}{2}$, $b_1 = \frac{7}{2} = \frac{f+1}{2}$. Thus, $G(\chi) = \sqrt{5}\cdots 5^2(a+bi)$, where $a^2 + b^2 = p = 5$. $n = (5^6-1)/28 = 558$, $s(n) = 10$, $t^{-1}(n) \equiv 2 \pmod 5$, $s(-n) = 14$, $t^{-1}(-n) = 3 \pmod 5$, $e_p = (\frac{5-1}{2}) \equiv 2 \pmod 5$. Since $\pm 2$ are the two roots of equation $x^2 + 1 = 0$ in $\mathbb{F}_5$, we can let $\phi(i) = 2 \in \mathbb{F}_p$ conditionally proper selection of $\mathfrak{P}$. By (3.10) and (3.11),
$$\begin{pmatrix} a \\ b \end{pmatrix} \equiv \begin{pmatrix} 1 & 2 \\ 1 & -2 \end{pmatrix}^{-1} \begin{pmatrix} 0 \\ (-1)^2 \cdot 2^2 \cdot \chi(-1) \end{pmatrix} \equiv \begin{pmatrix} -2 \\ 1 \end{pmatrix} \pmod p.$$

Therefore,
$$G(\chi_{\mathfrak{P}}) = -\sqrt{5} \cdot 5^2 (2-i), \quad G(\overline{\chi}_{\mathfrak{P}}) = -\sqrt{5} \cdot 5^2 (2+i).$$

When $p \equiv 3 \pmod 4$ and $a_1 = 2$, it belongs to Case F3(1). We take $p = 11$, then $<p>= \{1,11,9,15,25,23\}$, $g_1 <p>= \{3,5,13,17,19,27\}$ and $b_0 = b_1 = \frac{f}{2} = 3$. By Theorem 4.9(iii) in [8], we have that
$$G(\chi) = (-1)^{\frac{11+1}{4}} \cdot 11^3 = -11^3.$$

And by Theorem 3.1, we can verify that the sigh of $11^3$ is $(-1)^{\frac{f}{2}+1}/t(-n) \equiv (-1)^{3+1}/(-1) \equiv -1 \pmod 3$.

(3). Thirdly, let $l_1 = 11$, then $N = 4 \times 11 = 44$, $f = 10$, $h_1 = 1$.

When $p \equiv 1 \pmod 4$ and $a_1 = 1$, it belongs to Case F2(2). We take $p = 13$, then $<p>= \{1,13,37,41,5,21,9,29,25,17\}$ and $b_0 = \frac{9}{2} = \frac{f-1}{2}$, $b_1 = \frac{11}{2} = \frac{f+1}{2}$. Thus, $G(\chi) = \sqrt{13} \cdot 13^4(a+bi)$ where $a^2 + b^2 = p = 13$. And $n = (13^{10}-1)/44 = 3133147542$, $s(n) = 66$, $t^{-1}(n) \equiv 2 \pmod{13}$, $s(-n) = 54$, $t^{-1}(-n) = 7 \pmod{13}$, $e_p = (\frac{13-1}{2}) \equiv 5 \pmod{13}$. It is easy to know that $\pm 5$ are the two roots of equation $x^2 + 1 = 0$ in $\mathbb{F}_{13}$, and we can let $\phi(i) = 5 \in \mathbb{F}_p$, conditionally proper selection of $\mathfrak{P}$. By (3.10) and (3.11), we have
$$\begin{pmatrix} a \\ b \end{pmatrix} \equiv \begin{pmatrix} 1 & 5 \\ 1 & -5 \end{pmatrix}^{-1} \begin{pmatrix} (-1)^4 \cdot 5 \cdot 7 \cdot \chi(-1) \\ 0 \end{pmatrix} \equiv \begin{pmatrix} 2 \\ -3 \end{pmatrix} \pmod{13}.$$

Thus,
$$G(\chi_{\mathfrak{P}}) = \sqrt{13} \cdot 13^4 (2-3i), \quad G(\overline{\chi}_{\mathfrak{P}}) = \sqrt{13} \cdots 13^4 (2+3i).$$

When $p \equiv 3 \pmod 4$ and $a_1 = 2$, it's belong to Case F3(2). We take $p = 3$, then, $<p>= \{1,3,9,27,37,23,25,31,5,15\}$ and $h_1 = 1$. By Theorem 4.9(iii) in [8], we know that $a = \pm 5$, $b = \pm 1$, and $b_0(p-1) = s(-n) = 12$, $b_1(p-1) = s(n) = 8$, $t(n)^{-1} \equiv 1 \pmod 3$. It is easy to know that $\pm 1$ are the two roots of equation $x^2 + 11 = 0$ in $\mathbb{F}_3$, and we can let $\phi(\sqrt{-11}) = 1 \in \mathbb{F}_p$, conditionally proper selection of $\mathfrak{P}$. By (3.9) and (3.10), we have
$$\begin{pmatrix} a \\ b \end{pmatrix} \equiv \begin{pmatrix} 1 & 1 \\ 1 & -1 \end{pmatrix}^{-1} \begin{pmatrix} 0 \\ 1 \end{pmatrix} \equiv \begin{pmatrix} -1 \\ 1 \end{pmatrix} \pmod 3.$$

Thus, the Gauss sum of order 44 over $\mathbb{F}_3$ and its conjugation are:
$$G(\chi) = 3^4 \left( \frac{5+\sqrt{-11}}{2} \right), \quad G(\overline{\chi}) = 3^4 \left( \frac{5-\sqrt{-11}}{2} \right).$$

# 4 Sigh (root unit) ambiguities of Gauss sums in index 4 case

Recently, [2, 7, 3] have given explicit evaluation of Gauss sums in index 4 case. For the Gauss sums with odd order $N$ in the index 4 non-cyclic case, [7] solved the problem of sigh ambiguities by the congruences modulo certain odd factors $l$ of $N$. However, as we have said in previous text, when $l = 1, 3$, $\sqrt{-l} \in K$, which is the fixed subfield of $\sigma_p$, the integral ring $O_{\mathbb{Q}(\sqrt{-l})} \subset O_K$ has four or six unit roots, i.e., $< i >$ or $< \zeta_6 >$, they can not be distinguished by the congruences modulo $l$. And for the other situations of Gauss sums in index 4 case, we refer to [3] and [6].

As we known, [7] discussed many subcases and gave the results in two large tables. In this section, we give an unified theorem to solve the problem of sigh (unit root) ambiguities for all the subcases of Gauss sums with odd order $N$ in the index 4 non-cyclic case, by Stickelberger's congruence. And we give several examples, as complements for [7]. Particularly, we note that it occurs that $\sqrt{-3} \in K$ in some examples, and they can not be correctly calculated by the method in [7].

In this section, we always keep the following assumptions. Assume that

▼. $p$ is a prime number, $N \geqslant 2$, $(p(p-1), N) = 1$, and let the order of $p$ modulo $N$ is $f = \frac{\varphi(N)}{4}$. So, $N$ must be odd, $[(\mathbb{Z}/N\mathbb{Z})^* :< p >] = 4$, and the decomposition field $K$ of $p$ in $\mathbb{Q}(\zeta_N)$ is abelian quartic field.

▼. $q = p^f$, $\chi$ be a multiplicative character of order $N$ over $\mathbb{F}_q$.

▼. $-1 \notin < p > \subset (\mathbb{Z}/N\mathbb{Z})^*$, so $K$ is an imaginary field.

By Theorem.3.2 in [7], the Gauss sums with odd order $N$ in index 4 non-cyclic case have exactly the following three forms:

$$
\begin{aligned}
&(1). \ G(\chi) = \varepsilon p^{\frac{f}{2}}; \\
&(2). \ G(\chi) = 2^{-y} p^x (a + b\sqrt{-L}); \\
&(3). \ G(\chi) = 2^{-y} p^x (a + b\sqrt{-L})(a' + b'\sqrt{-L'}),
\end{aligned}
\tag{4.1}
$$

where $f = \frac{\varphi(N)}{4}$, $x, y \in \mathbb{Z}_+ \cup \{0\}$, $(a, b), (a', b') \in \mathbb{Z}^2$ are determined by certain quadratic positive equations, and $\sqrt{-L}, \sqrt{-L'} \in K$.

Similar as Section 3, for the Gauss sums as form (4.1)(1) and (4.1)(2), the sigh or unit root ambiguities can be solved by equations (3.4), (3.9) and (3.10) in finite field $\mathbb{F}_p$. We next consider the Gauss sums as form (4.1)(3)).

Suppose the quotient group $(\mathbb{Z}/N\mathbb{Z})^*/ < p > \cong \{1, \sigma, \tau, \sigma\tau\}$ and $\sigma(\sqrt{-L}) = -\sqrt{-L}, \sigma(\sqrt{-L'}) = \sqrt{-L'}$; $\tau(\sqrt{-L}) = \sqrt{-L}, \tau(\sqrt{-L'}) = -\sqrt{-L'}$. Since $\chi|_{\mathbb{F}_p} = 1$ is trivial character, applying Stickelberger's congruence, respectively, to $G(\chi)$, $\sigma(G(\chi))$, $\tau(G(\chi))$ and $\sigma\tau(G(\chi))$, we have that

$$
\begin{aligned}
G(\chi) &= (-1)^x \cdots 2^y \cdot \pi^{(p-1)x}(a + b\sqrt{-L})(a' + b'\sqrt{-L'}) \equiv \tfrac{-\pi^{s(-n)}}{t(-n)} \pmod{\mathfrak{P}^{s(-n)+1}}; \\
G(\chi^\sigma) &= (-1)^x \cdot 2^y \cdot \pi^{(p-1)x}(a - b\sqrt{-L})(a' + b'\sqrt{-L'}) \equiv \tfrac{-\pi^{s(-\sigma n)}}{t(-\sigma n)} \pmod{\mathfrak{P}^{s(-\sigma n)+1}}; \\
G(\chi^\tau) &= (-1)^x \cdot 2^y \cdot \pi^{(p-1)x}(a + b\sqrt{-L})(a' - b'\sqrt{-L'}) \equiv \tfrac{-\pi^{s(-\tau n)}}{t(-\tau n)} \pmod{\mathfrak{P}^{s(-\tau n)+1}}; \\
G(\chi^{\sigma\tau}) &= (-1)^x \cdot 2^y \cdot \pi^{(p-1)x}(a - b\sqrt{-L})(a' - b'\sqrt{-L'}) \equiv \tfrac{-\pi^{s(-\sigma\tau n)}}{t(-\sigma\tau n)} \pmod{\mathfrak{P}^{s(-\sigma\tau n)+1}}.
\end{aligned}
$$

By Stickelberger's theorem,

$$(p-1)b_0 = s(-n), \ (p-1)b_1 = s(-\sigma n), \ (p-1)b_2 = s(-\tau n), \ (p-1)b_3 = s(-\sigma\tau n),$$

and $x = \min\{b_0, b_1, b_2, b_3\}$. Let

$$
\begin{aligned}
T'_0 &= \tfrac{(-1)^{x+1}\cdots 2^y}{t(-n)}, \ T'_1 = \tfrac{(-1)^{x+1}\cdots 2^y}{t(-\sigma n)}, \ T'_2 = \tfrac{(-1)^{x+1}\cdots 2^y}{t(-\tau n)}, \ T'_3 = \tfrac{(-1)^{x+1}\cdots 2^y}{t(-\sigma\tau n)} \\
T_i &= \begin{cases} T'_i & \text{If } x = b_i; \\ 0 & \text{otherwise} \end{cases} \quad (i = 0, 1, 2, 3)
\end{aligned}
\tag{4.2}
$$

And we take $\pm\beta, \pm\beta'$ be the roots of equations $x^2 + L = 0$ and $x^2 + L' = 0$ in $\mathbb{F}_p$. Selecting proper $\mathfrak{P}$ such that $\phi(\sqrt{-L}) = \beta', \phi(\sqrt{-L'}) = \beta'$ under the action of Teichmüller map. Thus, we obtain the following equations in

$\mathbb{F}_p$,

$$\begin{pmatrix} aa' \\ a'b \\ ab' \\ bb' \end{pmatrix} = \begin{pmatrix} 1 & \beta & \beta' & \beta\beta' \\ 1 & -\beta & \beta' & -\beta\beta' \\ 1 & \beta & -\beta' & -\beta\beta' \\ 1 & -\beta & -\beta' & \beta\beta' \end{pmatrix}^{-1} \begin{pmatrix} T_0 \\ T_1 \\ T_2 \\ T_3 \end{pmatrix} \pmod{p} \qquad (4.3)$$

Summing up, we have the theorem.

**Theorem 4.1.** *For the Gauss sums with odd order $N$ in the index 4 non-cyclic case (see [3] for detailed classification), conditionally proper selection of $\mathfrak{P}$, the sigh (unit root) ambiguities can be solved by equations (3.4), (3.9) and (3.10), (4.2) and (4.3) in $\mathbb{F}_p$, according to the different forms as (4.1).*

Next, we give some examples, as supplements for [7].

**Example 2** (Case N-a)

(1). **Case(N-a1)**. Let $l_1 = 7$, $l_2 = 11$, then $N = 77$, $f = 15$, $\left(\frac{7}{11}\right) = -1$. We take $p = 37$, then $\mathrm{ord}_N(p) = 15$ and $h_2 = h(\mathbb{Q}(\sqrt{-11})) = 1$. Solving the quadratic positive equations, we get

$$4 \cdot 37 = a^2 + 11b^2, \; 37 \nmid ab \implies a = \pm 7, \; b = \pm 3.$$

By Stickelberger's theorem, $b_0(p-1) = s(-n) = 288$, $t(-n) \equiv 7 \pmod{37}$, $b_1(p-1) = s(n) = 252$, $t(n) \equiv 21 \pmod{37}$. And it's easy to know that $\pm 10$ are the two roots of $x^2 + 11 = 0$ in $\mathbb{F}_{37}$. Conditionally proper selection of $\mathfrak{P}$, we can take $\phi(\sqrt{-11}) = 10 \in \mathbb{F}_{37}$. By (3.9) and (3.10), we have that

$$\begin{pmatrix} a \\ b \end{pmatrix} \equiv \begin{pmatrix} 1 & 10 \\ 1 & -10 \end{pmatrix}^{-1} \begin{pmatrix} 0 \\ 2 \cdot 21^{-1} \end{pmatrix} \equiv \begin{pmatrix} 30 \\ 34 \end{pmatrix} \pmod{37}.$$

Therefore, $a = -7, b = -3$, i.e., the Gauss sum of order 77 over $\mathbb{F}_{37^{15}}$ and its conjugation are

$$G(\chi_{\mathfrak{P}}) = -\frac{1}{2} 37^7 \cdot (7 + 4\sqrt{-11}), \quad G(\overline{\chi}_{\mathfrak{P}}) = -\frac{1}{2} 37^7 \cdot (7 - 4\sqrt{-11}).$$

(2). **Case(N-a2)**. Let $l_1 = 11$, $l_2 = 5$, then $N = 55$, $f = 10$, $\left(\frac{l_1}{l_2}\right) = \left(\frac{11}{5}\right) = 1$ $h_{12} = 4$, $h_1 = 1$, $h_2 = 4$. We take $p = 59$, then

$$4 \cdot 59^2 = a^2 + 55b^2, \; 59 \nmid AB, \implies A = \pm 102, \; B = \pm 8.$$

By Stickelberger's theorem, $b_0(p-1) = s(-n) = 348$, $t(-n) \equiv 43 \pmod{59}$, $b_1(p-1) = s(n) = 232$, $t(n) \equiv 11 \pmod{59}$. And it's easy to know that $\pm 2$ are the two roots of $x^2 + 55 = 0$ in $\mathbb{F}_{59}$. Conditionally proper selection of $\mathfrak{P}$, we can take $\phi(\sqrt{-55}) = 2 \in \mathbb{F}_{59}$. By (3.9) and (3.10),

$$\begin{pmatrix} a \\ b \end{pmatrix} \equiv \begin{pmatrix} 1 & 2 \\ 1 & -2 \end{pmatrix}^{-1} \begin{pmatrix} 0 \\ -2 \cdot 11^{-1} \end{pmatrix} \equiv \begin{pmatrix} 16 \\ 51 \end{pmatrix} \pmod{59}.$$

Thus, $a = -102, b = -8$, i.e., the Gauss sum of order 55 over $\mathbb{F}_{59^{10}}$ and its conjugation are

$$G(\chi_{\mathfrak{P}}) = -59^8(102 + 8\sqrt{-55})/2, \; G(\chi_{\mathfrak{P}}) = -59^8(102 - 8\sqrt{-55})/2.$$

(3). **Case(N-a3)**. Let $l_1 = 7$, $l_2 = 5$, then $N = 35$, $f = 6$, $\left(\frac{l_1}{l_2}\right) = \left(\frac{7}{5}\right) = -1$ $h_{12} = 2$, $h_1 = 1$, $h_2 = 2$. We take $p = 79$, then

$$\begin{cases} 4 \cdot 79 = a^2 + 35b^2, \; 79 \nmid AB \\ 4 \cdot 79 = (a')^2 + 7(b')^2, \; 79 \nmid a'b' \end{cases} \implies \begin{array}{l} a = \pm 1, \; b = \pm 3; \\ a' = \pm 8, \; b' = \pm 6. \end{array}$$

Since $(\mathbb{Z}/N\mathbb{Z})^*/<p> \{1,2\} \times \{1,3\}$, and by Stickelberger Theorem, $s(-n) = 312$, $t(-n) \equiv 71 \pmod{79}$, $s(n) = 156$, $t(n) \equiv 69 \pmod{79}$, $s(-2n) = 234$, $t(n) \equiv 8 \pmod{79}$, $s(-3n) = s(2n) = 234$, $t(n) \equiv 10 \pmod{79}$.

And it's easy to know that $\pm 26$ and $\pm 25$ are respectively the roots of $x^2 + 35 = 0$ and $x^2 + 7 = 0$ in $\mathbb{F}_{79}$. Conditionally proper selection of $\mathfrak{P}$, we can take $\phi(\sqrt{-35}) = 26, \phi(\sqrt{-7}) = 25$. Since $\sigma_2$ fixes $\sqrt{-7}$ and $\sigma_3$ fixes $\sqrt{-35}$, by (4.2) and (4.3),

$$\begin{pmatrix} aa' \\ a'b \\ ab' \\ bb' \end{pmatrix} \equiv \begin{pmatrix} 8 \\ -24 \\ 6 \\ -18 \end{pmatrix} \pmod{79} \Rightarrow \begin{cases} a = 1 \\ b = -3 \\ a' = 8 \\ b' = 6 \end{cases}.$$

Thus, the Gauss sum of order 35 over $\mathbb{F}_{79^6}$ and its conjugation are

$$\begin{aligned} G(\chi_{\mathfrak{P}}) &= \tfrac{79^2}{2}(1 - 3\sqrt{-35})(4 + 3\sqrt{-7}); \\ G(\chi_{\mathfrak{P}}^2) &= \tfrac{79^2}{2}(1 + 3\sqrt{-35})(4 + 3\sqrt{-7}); \\ G(\chi_{\mathfrak{P}}^3) &= \tfrac{79^2}{2}(1 - 3\sqrt{-35})(4 - 3\sqrt{-7}); \\ G(\chi_{\mathfrak{P}}^{-1}) &= \tfrac{79^2}{2}(1 + 3\sqrt{-35})(4 - 3\sqrt{-7}). \end{aligned}$$

**Example 3**(Case(N-b)) Let $N = 231 = 11 \times 7 \times 3$, then $(11, 7, 3) \equiv (3, 3, 3) \pmod 4$, $f = 120/4 = 30$, $\left(\frac{7}{3}\right) = 1$, $\left(\frac{3}{11}\right) = 1 \left(\frac{7}{11}\right) = -1$. We know that there are 56 elements of order 30 modulo 231, which can be divided into seven groups. The primitive roots modulo 11,7 and 3 are respectively:

$$g_1 \in \{211, 127, 106, 85\} \pmod{11}, \ g_2 \in \{199, 166\} \pmod 7, \ g_3 = 155 \pmod 3.$$

(1). When $p \equiv g_1^2 g_2 g_3$, $p \in \{-79, -46, 5, 26, 38, 47, 59, 80\} \pmod{231}$. It's belong to Case(N-b2-$\alpha$2), and we take $p = 5$, $h_{123} = 12$. Solving the positive quadratic equations, we have

$$4 \cdot 5^{12} = a^2 + 231b^2, \ 5 \nmid ab, \Longrightarrow a = \pm 27886, \ b = \pm 928.$$

By Stickelberger's theorem, $s(-n) = 72$, $t(-n) \equiv 3 \equiv 2^{-1} \pmod 5$, $s(n) = 48$, $t(n) \equiv 2 \equiv 3^{-1} \pmod 5$. And $\pm 2$ are the two roots of equation $x^2 + 231 = 0$ in $\mathbb{F}_5$. Conditionally proper selection of $\mathfrak{P}$, we can take $\phi(\sqrt{-231}) = 2 \in \mathbb{F}_5$. By (3.9) and (3.10),

$$\begin{pmatrix} a \\ b \end{pmatrix} \equiv \begin{pmatrix} 1 & 2 \\ 1 & -2 \end{pmatrix}^{-1} \begin{pmatrix} 0 \\ 2 \cdot 2^{-1} \end{pmatrix} \equiv \begin{pmatrix} -1 \\ 3 \end{pmatrix} \pmod 5.$$

Thus, $a = -27886, b = 928$, i.e., the Gauss sum of order 231 over $\mathbb{F}_{5^{30}}$ and its conjugation are

$$G(\chi_{\mathfrak{P}}) = 5^9(-13943 + 464\sqrt{-231}), \ G(\overline{\chi}_{\mathfrak{P}}) = 5^9(-13943 - 464\sqrt{-231}).$$

(2). When $p \equiv g_1 g_2^2 g_3 \pmod{231}$, $p \equiv g_1 g_2^2 g_3 \in \{-103, -82, -31, 2, 74, 95, 107, 116\}$. It's belong to Case(N-b2-$\alpha$2). We take $p = 2$, then

$$4 \cdot 2^{12} = a^2 + 231 b^2, \ 2 \nmid ab \Longrightarrow a = \pm 103, b = \pm 5.$$

Now, we can not distinguish $\pm$ modulo 2, which means that we can not use Stickelberger's congruence and Theorem 4.1. So, we use the congruence modulo $l_1^{r_1}$, the method in [7]. $a \equiv 2p^{\frac{h_{123}}{2}} \equiv 2^6 \cdots 2 \equiv 2 \pmod 7$, so $a = -103$. Thus, the Gauss sum of order 231 over $\mathbb{F}_{2^{30}}$ and its conjugation are

$$\{G(\chi), G(\overline{\chi})\} = \{2^8(-103 \pm 5\sqrt{-231})\}.$$

(3). when $p \equiv g_1 g_2 g_3^2 \pmod{231}$, $p \equiv g_1 g_2 g_3^2 \in \{-86, -53, 19, 40, 52, 61, 73, 94\}$. It's also belong to Case(N-b2-$\alpha$2). But, we notice that $K = \mathbb{Q}(\sqrt{-3}, \sqrt{-231})$, and it implies that there must exist an ambiguity of $<\zeta_6>$ in the course of evaluation of Gauss sum. If we take $p = 19$, then,

$$4 \cdot 19^6 = a^2 + 231b^2, \ 19 \nmid AB, \Longrightarrow a = \pm 79884338, b = \pm 3271120.$$

So, $a \equiv \pm 3, b \equiv \pm 4$ (mod 19). By Stickelberger's theorem, we have $s(-n) = 324$, $t(-n) \equiv 15 \equiv 14^{-1}$ (mod 19), $s(n) = 216$, $t(n) \equiv 14 \equiv 15^{-1}$ (mod 19). Since $231 \equiv 3$ (mod 19), $\pm 4$ is the two roots of $x^2 + 231 = x^2 + 3 = 0$ in $\mathbb{F}_{19}$. Conditionally proper selection of $\mathfrak{P}$, we can take $\phi(\sqrt{-231}) = \phi(\sqrt{-3}) = 4 \in \mathbb{F}_{19}$. Then $\phi(\zeta_6) = \phi(\frac{1+\sqrt{-3}}{2}) = -7$ and $\phi(\zeta_6^{-1}) = 7$, $\phi(\zeta_3) = \phi(\zeta_6^2) = -8, \phi(\zeta_3^{-1}) = 8$. By (3.9) and (3.10),

$$\varepsilon \begin{pmatrix} a \\ b \end{pmatrix} \equiv \begin{pmatrix} 1 & 4 \\ 1 & -4 \end{pmatrix}^{-1} \begin{pmatrix} 0 \\ 2 \cdot 14 \end{pmatrix} \equiv \begin{pmatrix} -5 \\ 6 \end{pmatrix} \pmod{19}.$$

So $\phi(\varepsilon) = -8$, $a \equiv 3$, $b \equiv 4$ (mod 19). Thus, $\varepsilon = \zeta_3 = \frac{-1+\sqrt{-3}}{2}$, $a = -79884338$, $b = 3271120$, i.e., the Gauss sum of order 231 over $\mathbb{F}_{19^{30}}$ and conjugation are

$$G(\chi_\mathfrak{P}) = 19^9 \left(\frac{-1+\sqrt{-3}}{2}\right)(-39942169 + 1635560\sqrt{-231});$$
$$G(\overline{\chi}_\mathfrak{P}) = 19^9 \left(\frac{-1+\sqrt{-3}}{2}\right)(-39942169 - 1635560\sqrt{-231}).$$

(4). When $p \equiv g_1^2 g_2^2 g_3$ (mod 231), $p \in \{-94, -73, -61, -52, -40, -19, 53, 86\}$. It's belong to Case(N-b3-$\alpha$). Take $p = 53$, then $s(-n) = 780, t(-n) \equiv 1$ (mod 53). By (3.4),

$$\varepsilon \equiv (-1)^{30/2+1}/t(-n) \equiv 1 \pmod{53} \quad \Rightarrow \quad \varepsilon = 1.$$

Thus, the Gauss sum of order 231 over $\mathbb{F}_{53^{30}}$ is

$$G(\chi) = 53^{15}.$$

(5). When $p \equiv g_1^2 g_2 g_3^2$ (mod 231), $p \in \{-107, -95, -74, -2, 31, 82, 103, 115\}$. It's belong to Case(N-b3-$\alpha$). Take $p = 31 K = \mathbb{Q}(\sqrt{-11}, \sqrt{-3})$, then $s(-n) = 450, t(-n) \equiv 25 \equiv 5^{-1}$ (mod 31). By (3.4),

$$\varepsilon \equiv (-1)^{30/2+1}/t(-n) \equiv 5 \pmod{31}.$$

Since $11^2 \equiv -3$ (mod 31), conditionally proper selection of $\mathfrak{P}$, we can take $\phi(\sqrt{-3}) = 11 \in \mathbb{F}_{31}$, i.e., $\phi(\zeta_6) = 6, \phi(\zeta_3) = 5$. Thus, $\varepsilon = \zeta_3$, i.e., the Gauss sum of order 231 over $\mathbb{F}_{31^{30}}$ and its conjugation are

$$G(\chi_\mathfrak{P}) = \frac{-1+\sqrt{-3}}{2} 31^{15}, \quad G(\overline{\chi}_\mathfrak{P}) = \frac{-1-\sqrt{-3}}{2} 31^{15}.$$

(6). When $p \equiv g_1 g_2^2 g_3^2$ (mod 231), $p \in \{-80, -59, -47, -38, -26, -5, 46, 79\}$. It's belong to Case(N-b3-$\alpha$). Take $p = 79, K = \mathbb{Q}(\sqrt{-7}, \sqrt{-3})$, then $s(-n) = 1170, t(-n) \equiv 55 \equiv 23^{-1}$ (mod 79). By (3.4),

$$\varepsilon \equiv (-1)^{30/2+1}/t(-n) \equiv 23 \pmod{31}.$$

Since $32^2 \equiv -3$ (mod 79), we can take $\phi(\sqrt{-3}) = 32 \in \mathbb{F}_{79}$, conditionally proper selection of $\mathfrak{P}$, which means that $\phi(\zeta_6) = -23, \phi(\zeta_6^{-1}) = 23$. Thus, $\varepsilon = \zeta_6^{-1}$, i.e., the Gauss sum of order 231 over $\mathbb{F}_{79^{30}}$ and its conjugation are

$$G(\chi_\mathfrak{P}) = \frac{1-\sqrt{-3}}{2} 79^{15}, \quad G(\overline{\chi}_\mathfrak{P}) = \frac{1+\sqrt{-3}}{2} 79^{15}.$$

(7). The rest $(56 - 6 \times 8) = 8$ elements of order 30 modulo 231 are $\{-4, -16, -25, -37, -58, 17, 68, 101\}$. They are exactly generated by $p \equiv g_1 g_2 g_3$ (mod 231). But, we notice that $p^{15} \equiv -1$ (mod 231), i.e., $-1 \in <p>$. So, we can evaluate the Gauss sum by the formula of "pure" Gauss sums ([1, Thm11.6.3] or [8, Lemma2.2]) such as, if we take $p = 101$, then $101^{15} \equiv -1$ (mod 231), $t = 15, s = 1$. Thus, the Gauss sum of order 231 over $\mathbb{F}_{101^{30}}$ is

$$G(\chi) = \frac{101^{15}+1}{231} 101^{30/2} = 101^{15}.$$

**Acknowledgements** This work was supported by the National Science Fundamental of China (NSFC) with No. 10990011.

# References

[1] Berndt B.C., Evans R.J. and Williams K. S., *Gauss and Jacobi Sums*, New York: J.Wiley and Sons Company, 1997

[2] Feng K., Yang J. and Luo S., Gauss sum of Index 4: (1) cyclic case, *Acta Math. Sinica (Eng. Ser.)*, 2005, 21(6): 1425-1434

[3] Feng K. and Yang J., The evaluation of Gauss sums for characters of 2-power order in th index 4 case, accepted by Algebra Colloquium, 2009

[4] Ireland K. and Rosen M., A classical introduction to modern number theory, New York: Springer-Varlag, 1982

[5] Stickelberger L., Über eine Verallgemeinerung von der Kreistheilung, *Math. Ann.*, 1890, 37: 321-367

[6] Xia L., An application of Stickelberger's congruence,*Proceeding of the seventh international conference on information and management sciences*, 2008, 7: 330-336.

[7] Yang J., Luo S. and Feng K., Gauss sum of Index 4: (2) non-cyclic case, *Acta Math. Sinica (Eng. Ser.)*, 2006, 22(3): 833-844

[8] Yang J. and Xia L., Complete solving for explicit evaluation of Gauss sums in the index 2 case, submmitted to *Science in China Series A*, 2009. http://arxiv.org/abs/0911.5472